\documentstyle[12pt, page1]{article}
\newcommand{\sect}[1]{\section{#1}\setcounter{equation}{0}}

\font\mbn=msbm10 scaled \magstep1
\font\mbs=msbm7 scaled \magstep1
\font\mbss=msbm5 scaled \magstep1
\newfam\mbff
\textfont\mbff=\mbn
\scriptfont\mbff=\mbs
\scriptscriptfont\mbff=\mbss
\def\mbf{\fam\mbff}

\def\Q{{\mbf Q}}
\def\Z{{\mbf Z}}
\def\Co{{\mbf C}}

\newtheorem{Th}{Theorem}[section]
\newtheorem{Lm}[Th]{Lemma}

\newtheorem{D}[Th]{Definition}
\newtheorem{Proposition}[Th]{Proposition}
\newtheorem{R}[Th]{Remark}

\author{Alexander Brudnyi\\
Department of Mathematics and Statistics\\
University of Calgary\\
Calgary, Canada}
\title{Solvable Quotients of K\"{a}hler Groups}
\date{}
\begin{document}
\maketitle
\begin{abstract}
{We prove several results on the structure of 
solvable quotients of fundamental groups of compact K\"{a}hler 
manifolds ({\em K\"{a}hler groups}).}
\end{abstract}
\sect{\hspace*{-1em}. Introduction.} 
We first recall a definition from [AN].
\begin{D}\label{de1}
A solvable group $\Gamma$ has finite rank, if there is a decreasing sequence
$\Gamma=\Gamma_{0}\supset\Gamma_{1}\supset\dots\supset\Gamma_{m+1}=\{1\}$
of subgroups, each normal in its predecessor, such that 
$\Gamma_{i}/\Gamma_{i+1}$ is abelian and $\Q\otimes(\Gamma_{i}/\Gamma_{i+1})$
is finite dimensional for all $i$.
\end{D}
In what follows ${\cal F}_{r}$ denotes a free group with the number of 
generators $r\in\Z_{+}\cup\{\infty\}$. Our main result is
\begin{Th}\label{te1}
Let $M$ be a compact K\"{a}hler manifold. Assume that the 
fundamental group $\pi_{1}(M)$ is defined by the sequence
$$
\{1\}\longrightarrow F\longrightarrow\pi_{1}(M)\stackrel{p}{\longrightarrow} H
\longrightarrow\{1\}
$$
where $H$ is a solvable group of finite rank of the form
$$
\{0\}\longrightarrow A\longrightarrow H\longrightarrow B\longrightarrow
\{0\}
$$
with non-trivial abelian groups $A,B$ so that $\Q\otimes A\cong\Q^{m}$ and
$m\geq 1$. Assume also that $p^{-1}(A)\subset\pi_{1}(M)$
does not admit a surjective homomorphism onto ${\cal F}_{\infty}$.
Then all eigen-characters of the conjugate action of $B$ on the vector
space $\Q\otimes A$ are torsion.
\end{Th}
In Lemma \ref{le2} we will show
that the condition for $p^{-1}(A)$ holds
if $F$ does not admit a surjective homomorphism onto ${\cal F}_{\infty}$. 

Using Theorem \ref{te1} we prove a result on solvable quotients of
K\"{a}hler groups. 
\begin{Th}\label{te2}
Assume that a K\"{a}hler group $G$ is defined by the sequence
$$
\{1\}\longrightarrow F\longrightarrow G
\stackrel{q}{\longrightarrow} H\longrightarrow\{1\}
$$
where $F$ does not admit a surjective homomorphism onto ${\cal F}_{\infty}$
and $H$ is a solvable group of finite rank.
Then there exist normal subgroups $H_{1}\supset H_{2}$ of $H$ so that
\begin{description}
\item[{\rm (a)}]\ $H_{1}$ has finite index in $H$;
\item[{\rm (b)}]\ $H_{1}/H_{2}$ is nilpotent, and
\item[{\rm (c)}]\ $H_{2}$ is torsion.
\end{description}
\end{Th}
\begin{R}\label{re1}
{\rm (1) Clearly the conclusion of Theorem \ref{te1} is valid for $H$ being 
an extension of $\Z^{n}$ by $\Z^{m}$, $n,m\geq 1$, and $F$ being a finitely
generated group. Assume that not all eigen-characters of the 
action of $\Z^{n}$ on $\Q\otimes\Z^{m}$ are torsion. 
Then a semidirect product of
such $F$ and $H$ (where $F$ is a normal subgroup of this product) is not
a K\"{a}hler group.\\
(2) Let $G$ be a K\"{a}hler group. By $DG=D^{1}G$ we denote the derived 
subgroup of $G$, and set $D^{i}G=DD^{i-1}G$. Assume that
$H:=G/D^{n}G$, $n\geq 1$, is a solvable group of finite rank. Then it was 
proved in [AN, Th. 4.9] and [Ca, Th. 2.2] that $H$ satisfies conditions
(a)-(c) of Theorem \ref{te2}. It is a consequence of the fact that
$G$ does not admit a surjective homomorphism onto ${\cal F}_{r}$ with 
$2\leq r<\infty$.}
\end{R}
\sect{\hspace*{-1em}. Proof of Theorem 1.2.}
In what follows $T_{n}(\Co)\subset GL_{n}(\Co)$ denotes the Lie group of
upper triangular matrices. Let $T_{2}\subset T_{2}(\Co)$ be the Lie group of 
matrices of the form
\[
\left(
\begin{array}{cc}
a&b\\ 0&1\\
\end{array}
\right)\ \ \ \ \ (a\in\Co^{*},\ b\in\Co) ,
\]
$D_{2}\subset T_{2}$ and $N_{2}\subset T_{2}$ be the groups of diagonal 
and unipotent matrices.
Let $M$ be a compact K\"{a}hler manifold.
For a homomorphism $\rho\in Hom(\pi_{1}(M),T_{2})$ we let
$\rho_{a}\in Hom(\pi_{1}(M),\Co^{*})$ denote the upper diagonal character
of $\rho$. The main result used in our proofs is the following 
\begin{Proposition}\label{pro1}
Assume that $\pi_{1}(M)$ is defined by the sequence
$$
\{1\}\longrightarrow F\longrightarrow\pi_{1}(M)
\longrightarrow H\longrightarrow\{1\}
$$
where the normal subgroup $F$ does not admit a surjective
homomorphism onto ${\cal F}_{\infty}$.
Assume that $\rho\in Hom(\pi_{1}(M),T_{2})$ satisfies
$F\subset Ker(\rho_{a})$ but $F\not\subset Ker(\rho)$.
Then $\rho_{a}$ is a torsion character.
\end{Proposition}
{\bf Proof.} Given a character $\xi\in Hom(\pi_{1}(M),\Co^{*})$, let 
$\Co_{\xi}$ denote the associated $\pi_{1}(M)$-module. We define 
$\Sigma^{1}(M)$ to be the set of
characters $\xi$ such that $H^{1}(\pi_{1}(M),\Co_{\xi})$ is nonzero.
The structure of $\Sigma^{1}(M)$ was described in the consequent papers of
Beauville [Be], Simpson [S], Campana [Ca].\\
{\bf BSC Theorem.} {\em There is a finite number of surjective
holomorphic maps with connected fibres $f_{i}:M\longrightarrow C_{i}$
onto smooth compact complex curves of genus $\geq 1$
and torsion characters $\rho_{i},\xi_{j}\in Hom(\pi_{1}(M),\Co^{*})$ such 
that}
$$
\Sigma^{1}(M)=\bigcup_{i}\rho_{i}f_{i}^{*}Hom(\pi_{1}(C_{i}),\Co^{*})\cup
\bigcup_{j}\{\xi_{j}\}\ .
$$
Further, the group $N_{2}$ acts on $T_{2}$ by conjugation. Any two
homomorphisms from $Hom(\pi_{1}(M),T_{2})$ belonging to the orbit
of this action will be called {\em equivalent}.

Let $\rho\in Hom(\pi_{1}(M),T_{2})$ satisfy the conditions of Proposition
\ref{pro1}. Then it is well known that the class of equivalence of 
$\rho$ is uniquely defined by an element 
$c_{\rho}\in H^{1}(\pi_{1}(M),\Co_{\rho_{a}})$ (see e.g. [A, Prop. 2]).
In particular, if $c_{\rho}=0$, $\rho$ is equivalent to a representation into
$D_{2}$. In our case, $c_{\rho}\neq 0$ because $F\not\subset Ker(\rho)$.
Thus $\rho_{a}$ satisfies the conditions of BSC Theorem. If $\rho_{a}$
coincides with one of $\xi_{j}$ then it is torsion by the above theorem.
So assume that $\rho_{a}=\rho_{i}f_{i}^{*}\phi$ for some torsion 
character $\rho_{i}$ and $\phi\in Hom(\pi_{1}(C_{i}),\Co^{*})$. 
Let $K:=Ker(\rho_{i})\subset\pi_{1}(M)$ and $p: M_{1}\rightarrow M$ be
the Galois covering of $M$ corresponding to the finite abelian Galois group
$\pi_{1}(M)/K$. Let $M_{1}\stackrel{g}{\longrightarrow}C
\stackrel{h}{\longrightarrow}C_{i}$ be the Stein factorization
of $f_{i}\circ p$. Here $g$ is a morphism with connected fibres
onto a smooth curve $C$ and $h$ is a finite morphism.
Assume, to the contrary, that $\rho_{a}$ is not torsion. Then we prove
\begin{Lm}\label{le1}
$F\subset Ker(g_{*})$. 
\end{Lm}
{\bf Proof.}
Set $G:=(f_{i})_{*}(F)\subset\pi_{1}(C_{i})$. According to the assumptions
of Proposition \ref{pro1} we have either (a) $G$ is a subgroup of finite
index in $\pi_{1}(C_{i})$, or (b) $G$ is isomorphic to ${\cal F}_{r}$ with 
$r<\infty$. Let us consider (a). Since by definition
$F\subset Ker(\rho_{a})$ and $\rho_{i}$ is torsion, 
$\phi(G)$ is a finite abelian group. Then $G_{1}:=Ker\ \phi\cap G$ is
a subgroup of finite index in $G$. In particular, $G_{1}$ is a subgroup
of finite index in $\pi_{1}(C_{i})$ and so it is not free. But
by our assumption, $\phi$ is not torsion and so $Ker(\phi)\cong
{\cal F}_{\infty}$. Thus $G_{1}$ is also free as a subgroup of 
${\cal F}_{\infty}$. This shows that (a) is never happen.

Consider now (b). Using the fact that
$(f_{i})_{*}$ is a surjection, we conclude that $G$ is a normal 
subgroup of $\pi_{1}(C_{i})$. Let $S\rightarrow C_{i}$ be a regular 
covering corresponding to $Q:=\pi_{1}(C_{i})/G$. If $r\geq 2$ then the group 
$Iso(S)$ of isometries of $S$ (with respect to the hyperbolic metric) is 
finite and since $Q$ is infinite we have $r\leq 1$.
If $r=1$, any discrete subgroup of $Iso(S)$ is virtually cyclic and 
in particular does not act cocompactly on $S$. Thus $r=0$ which means
that $G=\{e\}$ and $F\subset Ker(f_{i})_{*}$. Then the assumption
of the proposition implies that $F\subset K=\pi_{1}(M_{1})$. 

Now consider $\tilde\rho:=\rho|_{\pi_{1}(M_{1})}$ with
the upper diagonal character 
$\tilde\rho_{a}:=\rho_{a}|_{\pi_{1}(M_{1})}$. Since 
$\pi_{1}(M_{1})\subset\pi_{1}(M)$ is a
subgroup of finite index, $\tilde\rho_{a}$ is also not 
torsion. Set $\tilde\phi:=h^{*}\phi$. Then 
$\tilde\rho_{a}=g^{*}\tilde\phi$.
Let $G_{1}:=g_{*}(F)\subset\pi_{1}(C)$. Then the same
argument for $G_{1}$ as above for $G$ (with $\tilde\rho_{a}$ and
$\tilde\phi$ instead of $\rho_{a}$ and $\phi$)
yields $F\subset Ker(g_{*})$.\ \ \ \ \ $\Box$

According to Lemma \ref{le1} and the assumptions of Proposition \ref{pro1}
we have that $\tilde\rho|_{Ker(g_{*})}$ is non-trivial and $\tilde\rho_{a}$ is 
the pullback of a character from $Hom(\pi_{1}(C),\Co^{*})$. Then
from [Br, Proposition 3.6] it follows that $\tilde\rho_{a}$ is torsion. 
Therefore $\rho_{a}$ is torsion, as well.
This contradiction proves the proposition.\ \ \ \ \ $\Box$

We are ready to prove Theorem \ref{te1}. According to the assumptions of
the theorem there is a homomorphism $i$ of $H$ into the Lie group $R$ of the 
form
$$
\{0\}\longrightarrow\Co^{m}\longrightarrow R\longrightarrow B
\longrightarrow\{0\}
$$
whose kernel is $Tor(A)$. Here we idenitify $\Co^{m}$ with $\Co\otimes A$.
Consider the action $s:B\longrightarrow GL_{m}(\Co)$ by conjugation. Since,
by the definition of $\pi_{1}(M)$, $B$ is a finitely generated abelian group,
$s=\oplus_{j=1}^{d}s_{j}$ where $s_{j}$ is equivalent to a nilpotent
representation $B\longrightarrow T_{m_{j}}(\Co)$ with a diagonal character
$\rho_{j}$. Here $\sum_{j=1}^{d}m_{j}=m$. From this decomposition it follows
that there is an invariant $B$-submodule $V_{j}\subset\Co^{m}$ of
$dim_{\Co}V_{j}=m-1$ such that $W_{j}=\Co^{m}/V_{j}$ is a one-dimensional
$B$-module and the action of $B$ on $W_{j}$ is defined as 
multiplication by the character $\rho_{j}$. By definition, $V_{j}$ is
a normal subgroup of $R$ and the quotient group $R_{j}=R/V_{j}$ is defined
by the sequence
$$
\{0\}\longrightarrow\Co\longrightarrow R_{j}\longrightarrow B
\longrightarrow\{0\}\ .
$$
Here the action of $B$ on $\Co$ is multiplication by the character
$\rho_{j}$. (As before the associated $B$-module is denoted by
$\Co_{\rho_{j}}$.)
Let us denote by $t_{j}$ the composite homomorphism
$\pi_{1}(M)\longrightarrow H\stackrel{i}{\longrightarrow} R\longrightarrow
R_{j}$. Further, the equivalence class of extensions of $B$ by $\Co$
isomorphic to $R_{j}$
is defined by an element $c_{j}\in H^{2}(B,\Co_{\rho_{j}})$. We 
assume that the character $\rho_{j}$ is non-trivial (for otherwise,
$\rho_{j}$ is clearly torsion). Then $H^{2}(B,\Co_{\rho_{j}})=0$
(for the proof see e.g. [AN, Lemma 4.2]). This shows that $R_{j}$ is
isomorphic to the semidirect product of $\Co$ and $B$, i.e.,
$R_{j}=\Co\times B$ with multiplication
$$
(v_{1},g_{1})\cdot (v_{2},g_{2})=
(v_{1}+\rho_{j}(g_{1})\cdot v_{2},g_{1}\cdot g_{2}),\ \ \
v_{1},v_{2}\in \Co,\ g_{1},g_{2}\in B\ .
$$ 
Let us determine a map $\phi_{j}$ of $R_{j}$ to $T_{2}$
by the formula
\[
\phi_{j}(v,g)=\left(
\begin{array}{cc}
\rho_{j}(g)&v\\ 0&1\\
\end{array}
\right)
\]
Obviously, $\phi_{j}$ is a correctly defined homomorphism with  upper 
diagonal character $\rho_{j}$. Hence 
$\phi_{j}\circ t_{j}:\pi_{1}(M)\longrightarrow T_{2}(\Co)$ is a homomorphism 
non-trivial on $p^{-1}(A)\subset\pi_{1}(M)$ by its definition.
Also $p^{-1}(A)\subset Ker(\rho_{j}\circ t_{j})$. Since by our assumptions 
$p^{-1}(A)$ does not admit a surjective homomorphism onto 
${\cal F}_{\infty}$, Proposition \ref{pro1} 
applied to $\phi_{j}\circ t_{j}$ implies that $\rho_{j}$ is torsion. This
completes the proof of the theorem.\ \ \ \ \ $\Box$

We prove now the following result.
\begin{Lm}\label{le2}
Assume that a group $G$ is defined by the sequence
$$
\{1\}\longrightarrow G_{1}\longrightarrow G\longrightarrow G_{2}
\longrightarrow\{1\}
$$
where $G_{1},G_{2}$ do not admit surjective homomorphisms onto 
${\cal F}_{\infty}$.  Then $G$ satisfies the similar property.
\end{Lm}
{\bf Proof of Lemma \ref{le2}.}  
Assume, to the contrary, that there is a surjective homomorphism
$\phi:G\longrightarrow {\cal F}_{\infty}$. Then $\tilde G_{1}:=\phi(G_{1})$ 
is a normal subgroup of ${\cal F}_{\infty}$ and 
$\tilde G_{2}:={\cal F}_{\infty}/\tilde G_{1}$ is a quotient of $G_{2}$. 
Now the assumption of the 
lemma implies that $\tilde G_{1}\cong {\cal F}_{r}$ with
$r<\infty$. Let $X$ be a complex hyperbolic surface with 
$\pi_{1}(X)={\cal F}_{\infty}$ and $S\longrightarrow X$ be the regular 
covering corresponding to $\tilde G_{2}$.  Assume first that $r\geq 1$.
Since $\pi_{1}(S)=\tilde G_{1}$, any subgroup of the group
$Iso(S)$ of isometries of $S$ (with respect to the hyperbolic metric)
is finitely generated. In 
particular $\tilde G_{2}$ is finitely generated (as well as 
$\tilde G_{1}$). This implies that ${\cal F}_{\infty}$ should be finitely 
generated which is wrong. Thus $r=0$ and $\tilde G_{2}={\cal F}_{\infty}$.
This contradicts to our assumption and shows that
there is no such $\phi$.
\ \ \ \ \ $\Box$
\sect{\hspace*{-1em}. Proof of Theorem 1.4.}
For a group $L$ set $L^{ab}:=L/DL$. We say that an $L$-module $V$ is
{\em quasi-unipotent} if there is a subgroup $L'\subseteq L$ of finite
index whose elements act unipotently on $V$. To prove the theorem we
will check the following condition from [AN, Lemma 4.8].
\begin{Lm}\label{le3}
Let $H'\subseteq H$ be a subgroup of finite index. Then $H'$ acts
quasi-unipotently on the finite-dimensional vector space 
$\Q\otimes (H'\cap DH)^{ab}$.
\end{Lm}
{\bf Proof.} We set $K:=H'/(H'\cap DH)$, $G':=q^{-1}(H')$, $S:=D(H'\cap DH)$,
and $S':=q^{-1}(S)$. Here $K$ is a finitely generated abelian 
group. Indeed, $H'$ is finitely generated as a subgroup of 
finite index of the finitely generated
group $H$ ($=$ the image of the finitely generated group $G$).
Thus $K$ is finitely generated as the image of $H'$. The fact that 
$K$ is abelian follows directly from the definition. Further,
since $H'$ is a subgroup of finite index in $H$, $G'$ is a subgroup of
finite index in $G$. In particular, it is K\"{a}hler. Moreover, we have
$$
\{1\}\longrightarrow S'\longrightarrow G'\longrightarrow L
\longrightarrow\{1\}
$$
where $L$ is a solvable group of finite rank defined by the sequence
$$
\{0\}\longrightarrow (H'\cap DH)^{ab}\longrightarrow L\longrightarrow K
\longrightarrow\{0\}\ .
$$
Now the statement of the lemma is equivalent to the fact that
all eigen-characters of the conjugate action of $K$ on 
$\Q\otimes (H'\cap DH)^{ab}$ are
torsion. To prove that it suffices to show that $S'$ does not admit a
surjective homomorphism onto ${\cal F}_{\infty}$, and then to apply 
Lemma \ref{le2} and Theorem \ref{te1}.

Note that $S'$ is defined by the sequence
$$
\{1\}\longrightarrow F\longrightarrow S'\stackrel{q}{\longrightarrow}
S\longrightarrow\{1\}
$$
where $S$ is a solvable group of finite rank. By our assumption
$F$ does not admit a surjective homomorphism onto ${\cal F}_{\infty}$. Thus 
by Lemma \ref{le2}, $S'$ satisfies the same property.\ \ \ \ \ $\Box$

Now Theorem \ref{te2} is the consequence of 
Lemma \ref{le3} and [AN, Lemma 4.8].\ \ \ \ \ $\Box$


\begin{thebibliography}{      }
\bibitem[A]{A}
M.F.Atiyah, Complex analytic connections in fibre bundles, Trans. Amer. Math.
Soc. {\bf 85} (1957), 181-207.
\bibitem[AN]{AN}
D.Arapura and M.Nori, Solvable fundamental groups of algebraic varieties
and K\"{a}hler manifolds, Comp. Math. {\bf 116} (1999), 173-188.
\bibitem[Be]{Be}
A.Beauville, Annulation du $H^{1}$ pour les fibr\'{e}s en droites plats,
Lecture Notes in Math., {\bf 1507}, Springer Verlag (1992), 1-15.
\bibitem[Br]{Br}
A.Brudnyi, Solvable matrix representations of K\"{a}hler groups, preprint
(2001), 27 pp.
\bibitem[Ca]{Ca}
F.Campana, Ensembles de Green-Lazarsfeld et quotients resolubles des
groupes de K\"{a}hler, J. Alg. Geom. {\bf 10} (2001), no. 4, 599-622.
\bibitem[S]{S}
C.T.Simpson, Subspaces of Moduli Spaces of rank one local systems, Ann.
Sc. ENS, {\bf 26} (1993), 361-401.
\end{thebibliography}
\end{document}